\documentclass[12pt]{article}
\usepackage{graphicx}
\usepackage{amsmath,amsthm,amssymb,enumerate}
\usepackage{euscript,mathrsfs}
\usepackage{color}
\usepackage{dsfont}
\usepackage{url}
\usepackage{esint}
\usepackage[left=2cm,right=2cm,top=3.5cm,bottom=3.5cm]{geometry}
\usepackage{color}
\usepackage[framemethod=tikz]{mdframed}
\allowdisplaybreaks

\usepackage{soul}

\catcode`\@=11 \@addtoreset{equation}{section}

\catcode`\@=12

\newtheorem{Theorem}{Theorem}[section]
\newtheorem{Proposition}[Theorem]{Proposition}
\newtheorem{Lemma}[Theorem]{Lemma}
\newtheorem{Corollary}[Theorem]{Corollary}

\theoremstyle{definition}
\newtheorem{Definition}[Theorem]{Definition}

\newtheorem{Remark}[Theorem]{Remark}

\newcommand{\bTheorem}[1]{
	\begin{Theorem} \label{T#1} }
	\newcommand{\eT}{\end{Theorem}}

\newcommand{\bProposition}[1]{
	\begin{Proposition} \label{P#1}}
	\newcommand{\eP}{\end{Proposition}}

\newcommand{\bLemma}[1]{
	\begin{Lemma} \label{L#1} }
	\newcommand{\eL}{\end{Lemma}}

\newcommand{\bCorollary}[1]{
	\begin{Corollary} \label{C#1} }
	\newcommand{\eC}{\end{Corollary}}

\newcommand{\bRemark}[1]{
	\begin{Remark} \label{R#1} }
	\newcommand{\eR}{\end{Remark}}

\newcommand{\bDefinition}[1]{
	\begin{Definition} \label{D#1} }
	\newcommand{\eD}{\end{Definition}}

\newcommand{\Del}{\Delta_x}

\newcommand{\tvuh}{\widetilde{\vu}_h} 
\newcommand{\tTheta}{\widetilde{\Theta}} 
\newcommand{\tvZ}{\widetilde{Z}}
\newcommand{\avintO}[1]{\fint_{\Omega} #1 \dx}

\newcommand{\Ds}{\mathbb{D}_x}

\newcommand{\data}{{\rm data}}

\newcommand{\bFormula}[1]{
	\begin{equation} \label{#1}}
	\newcommand{\eF}{\end{equation}}

\newcommand{\vun}{\vu_n}

\newcommand{\vuh}{\vu_h}

\newcommand{\Divh}{{\rm div}_h}
\newcommand{\Gradh}{\nabla_h}

\newcommand{\Ov}[1]{\overline{#1}}

\newcommand{\aleq}{\stackrel{<}{\sim}}

\newcommand{\vr}{\varrho}
\newcommand{\vre}{\vr_\ep}

\newcommand{\vte}{\vt_\ep}
\newcommand{\vue}{\vu_\ep}

\newcommand{\tvu}{{\wtilde u}}
\newcommand{\tvt}{\widehat \vt}

\newcommand{\vt}{\vartheta}
\newcommand{\vu}{\vc{u}}

\newcommand{\vc}[1]{{\bf #1}}

\newcommand{\Div}{{\rm div}_x}
\newcommand{\Grad}{\nabla_x}

\newcommand{\dx}{\,{\rm d} {x}}

\newcommand{\dt}{\,{\rm d} t }

\newcommand{\intO}[1]{\int_{\Omega} #1 \ \dx}

\newcommand{\intOh}[1]{\int_{\Omega_h} #1 \ \dx}

\newcommand{\D}{{\rm d}}

\newcommand{\ep}{\varepsilon}

\newcommand{\vtB}{\vt_B}

\newcommand{\br}{ \nonumber \\ }

\def\softd{{\leavevmode\setbox1=\hbox{d}%
		\hbox to 1.05\wd1{d\kern-0.4ex{\char039}\hss}}}
\definecolor{Cgrey}{rgb}{0.85,0.85,0.85}
\definecolor{Cblue}{rgb}{0.50,0.85,0.85}
\definecolor{Cred}{rgb}{1,0,0}
\definecolor{fancy}{rgb}{0.10,0.85,0.10}
\definecolor{amaranth}{rgb}{0.9, 0.17, 0.31}

\newcommand\Cbox[2]{%
	\newbox\contentbox%
	\newbox\bkgdbox%
	\setbox\contentbox\hbox to \hsize{%
		\vtop{
			\kern\columnsep
			\hbox to \hsize{%
				\kern\columnsep%
				\advance\hsize by -2\columnsep%
				\setlength{\textwidth}{\hsize}%
				\vbox{
					\parskip=\baselineskip
					\parindent=0bp
					#2
				}%
				\kern\columnsep%
			}%
			\kern\columnsep%
		}%
	}%
	\setbox\bkgdbox\vbox{
		\color{#1}
		\hrule width  \wd\contentbox %
		height \ht\contentbox %
		depth  \dp\contentbox
		\color{black}
	}%
	\wd\bkgdbox=0bp%
	\vbox{\hbox to \hsize{\box\bkgdbox\box\contentbox}}%
	\vskip\baselineskip%
}

\mdfdefinestyle{MyFrame}{%
	linecolor=black,
	outerlinewidth=1pt,
	roundcorner=5pt,
	innertopmargin=\baselineskip,
	innerbottommargin=\baselineskip,
	innerrightmargin=10pt,
	innerleftmargin=10pt,
	backgroundcolor=white!20!white}

\newcommand{\wtilde}{\widetilde}

\allowdisplaybreaks

\begin{document}


\title{\bf On a system of equations arising in meteorology: Well-posedness and data assimilation }

\author{Eduard Feireisl
	\thanks{The work of E.F. was partially supported by the
		Czech Sciences Foundation (GA\v CR), Grant Agreement
		24--11034S. The Institute of Mathematics of the Academy of Sciences of
		the Czech Republic is supported by RVO:67985840. This work was partially supported by the Thematic Research Programme, University of Warsaw, Excellence Initiative Research University.} \and Piotr Gwiazda \thanks{The work of P. G.  and A. \'S.-G.  was partially supported by  National Science Centre
		(Poland),  agreement no   2024/54/A/ST1/00159.}
	\and Agnieszka \'Swierczewska-Gwiazda 
	}

\date{}

\maketitle

\medskip

\centerline{Institute of Mathematics of the Academy of Sciences of the Czech Republic}

\centerline{\v Zitn\' a 25, CZ-115 67 Praha 1, Czech Republic}

\centerline{feireisl@math.cas.cz}

\medskip

\centerline{Institute of Mathematics of Polish Academy of Sciences}
\centerline{\'Sniadeckich 8, 00-956 Warsaw, Poland}
\centerline{and}
\centerline{Interdisciplinary Centre for Mathematical and Computational Modelling, University of Warsaw}
\centerline{Pawi\'skiego 5a, 02-106 Warsaw, Poland}
\centerline{pgwiazda@mimuw.edu.pl}

\medskip

\centerline{Institute of Applied Mathematics and Mechanics, University of Warsaw}

\centerline{Banacha 2, 02-097 Warsaw, Poland}
\centerline{aswiercz@mimuw.edu.pl}

\begin{abstract}
	

Data assimilation plays a crucial role in modern weather prediction, providing a systematic way to incorporate observational data into complex dynamical models.  The paper addresses continuous data assimilation for a model arising as a singular limit of the three-dimensional compressible Navier-Stokes-Fourier system with rotation driven by temperature gradient. The limit system preserves the essential physical mechanisms of the original model,
while exhibiting a reduced, effectively two-and-a-half-dimensional structure. This simplified framework allows for a rigorous analytical study of the data assimilation process while maintaining a direct physical connection to the full compressible model. We establish well posedness of global-in-time solutions and a compact trajectory attractor, followed by the stability and convergence results for the nudging scheme applied to the limiting system. Finally, we demonstrate how these results can be combined with a relative entropy argument to extend the assimilation framework to the full three-dimensional compressible setting, thereby establishing a rigorous connection between the reduced and physically complete models.

\end{abstract}


{\small

\noindent

\medbreak

}

\section{Introduction}
\label{i}

Data assimilation is an essential approach used to improve the accuracy of numerical simulations, particularly in complex and nonlinear systems such as equations modelling the atmosphere. Reliable predictions depend strongly on the precision of the initial conditions, yet in practice these conditions can only be measured discretely. Weather observations, for example, are collected from satellites, ground stations, LiDAR and radar, but these measurements provide only partial snapshots of the true state of the atmosphere. Unfortunately, the spatial resolution of these measurements is insufficient to accurately reconstruct the initial conditions on the computational grid.
Data assimilation provides a systematic framework to nudge model trajectories toward observational measurements. Since observations are collected over a time interval rather than at a single moment, they help to compensate for the missing spatial information and allow for the reconstruction of an initial condition that leads to a physically meaningful model evolution.

A central mathematical idea behind data assimilation for dissipative systems is that only a limited number of parameters or modes may be sufficient to determine the long-term behavior of a system. 

{Properties such as existence and uniqueness of solutions, 
are of fundamental importance in data assimilation, providing  that the model dynamics are well posed, so that observations can be consistently incorporated without ambiguity. For a broad class of models, the long-term behavior of solutions can be characterized using only a finite number of degrees of freedom, the first rigorous proof for two-dimensional Navier-Stokes is due to \cite{FoiasProdi}. The result demonstrates that, once the first $N$
Fourier modes of two Navier-Stokes solutions agree for sufficiently large 
$N$, the entire solutions necessarily converge to the same limit as time tends to infinity. 
To specify what {\it sufficiently large $N$} means, the paper \cite{Foias} provides an explicit criterion, which was later refined in \cite{JonesTiti}. Interestingly, \cite{FoiasTemam} demonstrated that an analogous result holds when the degrees of freedom are taken to be quantities other than Fourier modes; in their work, the values of the solution at selected nodes played this role. Subsequent studies \cite{FoaisTiti, JonesTiti2} replaced nodal values with averages over finite volumes. Moreover, \cite{JonesTiti3} established that, for the two-dimensional Navier-Stokes equations, a broad family of determining sets of degrees of freedom exists,  including those arising naturally in finite element discretizations.}


Various data assimilation methods have been developed to achieve this goal, including variational approaches such as the four-dimensional variational data assimilation (4D-Var) method and ensemble-based techniques like the Ensemble Kalman Filter (EnKF). Another important class of methods is based on the Azouani, Olson and Titi (AOT)  framework, cf. \cite{AzOlTi}, which introduces a continuous data assimilation scheme formulated directly at the level of the evolution equations. 
Azouani, Olson, and Titi  showed that imposing suitable constraints on the finite-dimensional resolution of the measurement data is sufficient to guarantee that the solution reconstructed by their algorithm converges, over time, to the unknown reference state. 
Instead of incorporating the observations directly into the nonlinear dynamics of the model, a feedback control mechanism is introduced, which steers the system toward the reference solution consistent with the data. This approach is motivated by the practical difficulty of embedding raw measurements into the governing equations. 
 To explain the idea in more detail
let $u(t)$ denote the reference solution of the system
\begin{equation}
    \frac{du}{dt} = F(u), \qquad u(0)=u_0,
\end{equation}
where the initial condition $u_0$ is unknown. Observations are available at coarse
spatial resolution through $I_h(u(t))$, for $t \in [0,T]$. The objective is to construct
an approximate solution $v(t)$, emerging from arbitrary initial data $v(0)=v_0$, such that
$v(t)$ converges to $u(t)$ as $t$ increases, thereby recovering a suitable initial
state for forecasting beyond $T$.

The algorithm is given by
\begin{equation}
    \frac{dv}{dt} = F(v) - \mu I_h(v) + \mu I_h(u), \qquad v(0)=v_0,
\end{equation}
where $\mu>0$ is a relaxation parameter. The feedback term
$\mu(I_h(u)-I_h(v))$ forces the coarse spatial scales of $v$ toward those of the
observations, ensuring convergence to the reference solution provided that $\mu$ and
the spatial resolution $h$ are chosen appropriately.

{The analytical expectations outlined above have motivated a focus on rigorous results for the two-dimensional Navier-Stokes system, as in \cite{AzOlTi}.  The later result \cite{FaJoTi} showed that, for 2D Rayleigh-B\'enard convection, measuring and assimilating (nudging) a single velocity component suffices to recover the exact reference solutions for both the full velocity field and the temperature. In other words, there is no need to measure or nudge the temperature, see also \cite{FarhatLunasinTiti2017} . Moreover, in \cite{AltafEtAl2017} the authors showed that, for the same model, nudging the temperature alone does not suffice to recover the velocity field.  Nevertheless, \cite{FarhatLunasinTiti2016} demonstrates that, in the case of the planetary geostrophic model, observing and nudging only the temperature already suffices to reconstruct both the complete velocity field and the temperature. This highlights that the appropriate choice of state variables to observe and nudge is inherently model dependent.
Follow-up work \cite{CaLaTi} explored sensitivity analyses and established super-exponential convergence rates for 2D Navier-Stokes data assimilation. Meanwhile, \cite{LunasinTiti2017} extended the assimilation framework to settings with partial observability and feedback control of determining parameters. The works  \cite{BiswasBrownMartinez2022} -- \cite{BiswasPrice2021} further broadened the theory by introducing mesh-free interpolant observables and local-observable schemes, as well as by extending nudging approaches to both 2D and 3D incompressible Navier-Stokes systems not including temperature effects. Note that the 3D case considered in \cite{BiswasPrice2021}  
the assimilation process is conditional under some additional requirement of the observables, e.g. this may fail to hold for observables generated by the Leray Hopf solutions.}

  In the present work, we examine the rotating Oberbeck-Boussinesq system and demonstrate how continuous data assimilation, applied through the relative entropy framework, can be used to provide insights into the full three-dimensional compressible model.
Our attention is directed to  a system of equations arising as a singular limit of three-dimensional fast rotating compressible viscous fluid driven by temperature gradient. The considered system is particularly interesting in the context of continuous data assimilation, as it originates from the limit of a three-dimensional model and, due to the effect of the centrifugal force, reduces to a two-and-a-half dimensional formulation, see \cite{FanFei2025} for rigorous estimates based on the relative entropy. This feature is especially relevant for data assimilation techniques, since allows to provide existence of unique solutions and global attractor. 
What is however important, the effective reduction in dimensionality is not imposed artificially but rather emerges naturally from the underlying dynamics due to fast rotation of the system. Moreover, the  Coriolis force is indeed present in models actually used in weather forecasting, for which data assimilation plays a crucial role, see for e.g. the Unified Model distributed by the UK's national meteorological service Met Office cf.~\cite{White}, as well as its local versions. e.g  Polish Meteo service (https://www.meteo.pl/) of the Interdisciplinary Centre for Mathematical and Computational Modelling of the University of Warsaw. Even though the Coriolis force appears solely in the primitive system considered here and is absent from the target system, it has a significant influence on the asymptotic behavior, to the extent that it effectively removes all vertical motion from the limiting dynamics.

The Oberbeck-Boussinesq approximation provides a physically consistent yet mathematically  and numerically tractable framework for modeling rotating stratified flows. It captures the essential features of atmospheric dynamics, while striking a balance between physical realism and mathematical tractability, avoiding the full complexity of the compressible Navier-Stokes system. In the present work, we focus on the so-called rotating Oberbeck-Boussinesq system, which is incompressible and accounts for the horizontal components of the velocity field together with three components of the temperature. For this reason we will often call it {\it two-and-a-half-dimensional model}.

  The system arises as a singular limit of compressible three-dimensional Navier-Stokes-Fourier system. One could, however, argue that the most realistic model for describing atmospheric behavior is given by the three-dimensional compressible equations, and therefore data assimilation should ideally be performed at that level. Interestingly, there is no need to adopt this point of view, as the combination of the relative entropy method and the results presented here for the rotating Oberbeck-Boussinesq system shows that, for well-prepared initial data, it is possible to assimilate data for the three-dimensional model fully with the aid of the considered two-and-a-half-dimensional one while only controlling the assimilation error.   This connection is discussed in detail in Section 5, following the detailed analysis (Theorems~\ref{eP1} and~\ref{LP3}) and rigorous data assimilation procedure (Theorem~\ref{tt1}) carried out for the simpler system in Sections 2-4. An important advantage of avoiding data assimilation directly in the three-dimensional regime is the significant reduction in computational cost, given that data assimilation alone may account for a significant  part of the total computational expense. This reduction may open a possibility to substantially lower the energy consumption and, consequently, the overall cost of weather prediction.

\smallskip

\subsection{Primitive system}

Here and herafter, the spatial variables are denoted as 
\[
x = (x_h, x_3),\ x_h \in R^2,\ x_3 \in R.
\]

Consider a viscous, compressible, heat conducting and fluid confined to rotating cylinder 
\begin{equation} \label{s2}
	\Omega = B_r \times (0,1),\ B(r) = \left\{ x_h \in R^2 \ \Big| \ |x_h| < r \right\}
\end{equation}
The density $\vr$, the absolute temperature $\vt$, and the velocity $\vu$ of the fluid are
governed by the  Navier--Stokes--Fourier system
written in the rotation frame:  
\begin{align} 
\partial_t \vr + \Div (\vr \vu) &= 0, \br
\partial_t (\vr \vu) + \Div (\vr \vu \otimes \vu) + 
\frac{1}{\sqrt{\ep}} \vr \vu \times e_3 + \frac{1}{\ep^2} \Grad p(\vr, \vt) &= \Div \mathbb{S}(\vt, \Ds \vu) + \frac{1}{\ep} \vr \Grad F, \br
\partial_t (\vr e(\vr, \vt)) + \Div (\vr e(\vr, \vt) \vu) + 
\Div \vc{q}(\vt, \Grad \vt)) &= \ep^2 \mathbb{S}(\vt, \Ds \vu): \Ds\vu - 
p(\vr, \vt) \Div \vu,
\label{s1}
\end{align}
with the boundary conditions
\begin{align} 
\vu &= 0 \ \mbox{for}\ |x_h| = r,\ x_3 \in (0,1), \br
\vu \cdot \vc{n} &= 0,\ 
(\mathbb{S}(\vt, \Ds \vu) \cdot \vc{n}) \times \vc{n} = 0 
\ \mbox{for}\ x_h \in B(r),\ x_3 = 0,1, \br	
\vt|_{\partial \Omega} &= \Ov{\vt} + \ep \vtB,\ \Ov{\vt} > 0,
	\label{s3}
\end{align} 
and with the prescribed total mass 
\begin{equation} \label{s4}
\intO{ \vr } = \Ov{\vr} |\Omega| .
\end{equation}
The forcing potential $F$ takes the form 
\begin{equation} \label{i5}
F = \vc{g} \cdot x + |x_h|^2,
\end{equation}	
where the vector $\vc{g}$ determines the direction and magnitude of the gravitational field, while $|x_h|^2$ represent a potential of the centrifugal force, see \cite{FanFei2025} for details concerning the equation of state as well as other relevant constitutive 
equations.	

The fluid motion is driven by a
simultaneous action of three forces: the gravitational force, the buoyancy force induced by the temperature gradient, and 
the rotation represented by the Coriolis and centrifugal forces. Models of this type arise frequently in meteorology, geophysics or astrophysics. 
The system of equations \eqref{s1} is scaled by a small parameter $\ep$: 
\begin{itemize}
\item The origin of the thermal convection flow is the \emph{compressibility} of the fluid. 
The Mach number is small proportional to $\ep$, meaning the fluid is nearly incompressible.
In the incompressible limit, where the pressure 
becomes constant, the density variation may be replaced by temperature variation with the opposite sign yielding the celebrated Boussinesq relation.  

\item For the limit system to feel this change, the Mach and Froude characteristic numbers must be properly scaled, cf. \cite{FeNo6A}. 

\item If the limit system is influenced by the centrifugal force, the same scaling must be applied. However, the scaled centrifugal force imposes imperatively 
smallnes of the Rossby number in the Coriolis force. 

\item In the regime of small Rossby number, the Taylor--Proudman theorem applies, enforcing the motion of the fluid to become purely horizontal.
This is in sharp contrast with the conventional and frequently used Oberbeck--Boussinesq approximation advocated by Ecke and Shishkina \cite{EckeS} and the references cited therein.
\end{itemize}

\subsection{Singular limit and the target system}

The problem considered in the present paper arises as a singular limit of solutions $(\vre, \vte, \vue)$ for $\ep \to 0$. As shown in \cite{FanFei2025}, 
\begin{align} 
	\vue &\to [\vuh, 0] ,\ \vuh = \vuh(t,x_h) \br
\frac{\vte - \Ov{\vt}}{\ep} - \lambda(\Ov{\vr}, \Ov{\vt}) 
\avintO{ \frac{\vte - \Ov{\vt}}{\ep}  } &\to \theta,\ \lambda (\Ov{\vr}, \Ov{\vt}) = 1 - \frac{c_v (\Ov{\vr}, \Ov{\vt})}{c_p(\Ov{\vr}, \Ov{\vt}) } \in (0,1)
\label{s5}	
	\end{align}
where the limit velocity $\vuh = [u^1_h, u^2_h,0]$ is purely horizontal and satisfies the incompressible Navier--Stokes system 
\begin{align} 
	\Divh \vuh &= 0, \br
\Ov{\vr} \partial_t \vu_h + \Ov{\vr} \Divh( \vuh \otimes \vuh) + \Gradh \Pi &=
\mu \Delta_h \vuh - \left< \theta \right > \Gradh F,\ \mu > 0,
\label{i3} 
\end{align}
in $(0,T) \times B_r$, with the no-slip boundary conditions  
\begin{equation} \label{i7}
	\vuh|_{\partial B_r} = 0.  
\end{equation}
Here, we have denoted
\begin{equation} \label{i4}
< \theta >(t, x_h) = \int_0^1 \theta(t, x_h, x_3) \ \D x_3. 
\end{equation}

The temperature deviation $\theta$ solves a perturbed heat equation 
\begin{equation} \label{i6} 
\partial_t \theta + \Div (\theta \vuh) = \kappa \Del \theta +
a \Div (F \vuh), \ \kappa > 0,\ a \in R, 
\end{equation}
in $(0,T) \times \Omega$,  with a nonlocal boundary condition identified in \cite{BelFeiOsch}:
\begin{equation} \label{i8}
\theta|_{\partial \Omega} = \vtB - \alpha \frac{1}{|\Omega|} \intO{ 
	\theta },\ \alpha = \frac{c_p(\Ov{\vr}, \Ov{\vt} )}{c_v(\Ov{\vr}, \Ov{\vt} )} - 1.
\end{equation}
We suppose that the boundary temperature $\vtB$ is a restriction of a smooth function 
on $\partial \Omega$, 
\begin{equation} \label{i9}
	\vtB \in C^2(R^3). 
\end{equation}
The constants $c_p$ and $c_v$ stand for the specific heat at constant pressure and volume, respectively, 
evaluated at a thermal equilibrium state $(\Ov{\vr}, \Ov{\vt} )$. 

We call the problem \eqref{i3}--\eqref{i8} \emph{rotating Oberbeck--Boussinesq (rOB) system}. 

\subsection{Mathematics of rotating Oberbeck--Boussinesq  system }

Our main goal is to develop a mathematical theory for the rOB system.

\begin{itemize}
\item In section \ref{e}, we show rOB system is well posed in the class of strong solutions for arbitrarily large initial/boundary data.

\item In Section \ref{L}, we discuss the long--time behavior of solutions. We show the system is dissipative in the sense of Levinson 
(cf. Haraux \cite{Hara1}), meaning it admits a bounded absorbing set and the corresponding trajectory attractor. 

\item In Section \ref{D}, we propose a continuous data assimilation method in the spirit of Azouani, Olson, and Titi \cite{AzOlTi}.
\item {Finally, in Section \ref{S5} we show that, for well-prepared initial data, the three-dimensional model can be fully assimilated using the rOB system, with only the assimilation error being controlled, however on finite time interval.}
\end{itemize}


\section{Well-posedness}
\label{e}

Our first objective is to establish existence of global--in--time strong solutions to the rOB system. Note the problem exhibits 
certain similarity with the ``genuine'' Oberbeck--Boussinesq system studied in \cite{AbbFei22}. We consider a general situation,
where the fluid is confined to a cylindrical domain.
\begin{equation} \label{i1}
\Omega = \Omega_h \times (0,1), 
\end{equation}
where 
\begin{equation} \label{i2}
\Omega_h \subset R^2 \ \mbox{is bounded and smooth of class}\ C^2.	
\end{equation}	

\subsection{Transformation} 
\label{t}

Without loss of generality, we may assume 
\begin{equation} \label{t1}
\Ov{\vr} = 1, \ \intO{ F } = 0. 
\end{equation}
Introducing 
\[
\Theta = \theta -  a F + 2 a x_3^2 - C,\ C = 2a \left< x_3^2 \right>,
\]
we may rewrite the system of equations \eqref{i3}, \eqref{i6} 
in the form 
\begin{align} 
	\Divh \vuh &= 0, \br
	\partial_t \vu_h + \Divh( \vuh \otimes \vuh) + \Gradh \Pi &=
	\mu \Delta_h \vuh - \left< \Theta \right > \Gradh F
	\label{t2} 
\end{align}
in $(0,T) \times \Omega_h$,
and 
\begin{equation} \label{t3} 
	\partial_t \Theta + \Div (\Theta \vuh) = \kappa \Del \Theta 
\end{equation}
in $(0,T) \times \Omega$, 
with the boundary conditions 
\begin{equation} \label{t4}
\vuh|_{\partial \Omega_h} = 0,\ 
\Theta|_{\partial \Omega} = \vtB - \alpha \avintO{ \Theta }, 
\end{equation}
where we have replaced the boundary data
\[
\vtB \approx \vtB - a F + 2 a x_3^2 - C.
\]

\subsection{Global existence}
\label{Ex}

Problem \eqref{t2}--\eqref{t4}, supplemented with the initial data 
\begin{equation} \label{e1}
\vuh(0, \cdot) = \vu_{h,0},\ \Theta(0,\cdot) = \Theta_0 
\end{equation}	
is well posed in the class of strong solutions. 

\begin{mdframed}[style=MyFrame]
 
\begin{Theorem}[{\bf Well posedness}] \label{eP1}
Let $\Omega_h \subset R^2$ be bounded domain of class $C^2$. Suppose the boundary temperature is a restriction of a smooth function $\vtB \in C^2(R^3)$.  
Let the initial data belong to the class
\begin{equation} \label{e2}
\vu_{0,h} \in W^{2 - \frac{2}{p} ,p}(\Omega_h),\ \Theta_0 \in W^{2 - \frac{2}{p} ,p} \cap W^{2,2}(\Omega) \ \mbox{for some}\ p \geq 2, 
\end{equation}
and satisfy the compatibility conditions 
\[
\vu_{0,h}|_{\partial \Omega_h} = 0, \ 
\Theta_0|_{\partial \Omega} = \vtB - \alpha \avintO{ \Theta_0 }. 	
\]	

Then there exists a global-in-time solution $\vuh$, $\Theta$ of 
problem \eqref{t2}, \eqref{t3}, supplemented with the boundary conditions \eqref{t4}, and the initial conditions \eqref{e1}, unique in the maximal regularity class 
\begin{align} 
\vuh &\in L^p(0,T; W^{2,p}(\Omega_h; R^2)),\ 
\partial_t \vuh \in L^p(0,T; L^p(\Omega_h; R^2)), \br
\Theta 	&\in L^p(0,T; W^{2,p}(\Omega)) \cap L^\infty(0,T; W^{2,2}(\Omega)),\ 
\partial_t \Theta \in L^p(0,T; L^p(\Omega)) \cap L^\infty(0,T; L^2(\Omega)),
\label{e3}
\end{align}	
$T > 0$.	
\end{Theorem}

\end{mdframed}

The rest of this section is devoted to the proof of Theorem \ref{eP1}. We focus only on {\it a priori} bounds, the solution can be then constructed by the method described in detail in \cite{AbbFei22}.

\subsection{\emph{A priori} bounds}
\label{a}

We derive suitable {\it a priori} estimates for problem \eqref{t2}--\eqref{t4}. 

\subsubsection{Harmonic extension of boundary data}
\label{h}

First, we introduce the (unique) harmonic extension of the boundary 
data
\begin{equation} \label{a1}
	\tvt_B|_{\partial \Omega} = \vtB,\ 
	\Del \tvt_B = 0 \ \mbox{in}\ \Omega.
\end{equation}	
As $\Omega$ is merely Lipschitz, regularity of $\tvt_B$ is not {\it a priori} guaranteed. However, we can write 
\[
\tvt_B = \vtB + \eta, 
\]
where 
\begin{equation} \label{a2}
\Del \eta = - \Del \vtB \ \mbox{in}\ \Omega,\ \eta|_{\partial \Omega} = 0.
\end{equation}
Extending $- \Del \vtB$ as an odd function periodic in $x_3$, we may look for $\eta$ as an odd $x_3$-periodic solution of the problem 
\begin{equation} \label{a3}
\Del \eta = - \Del \vtB \ \mbox{in}\ \Omega_h \times [-1,1]_{\{ -1, 1\} }.	
\end{equation}
As $\vtB$ is smooth, the extension of $\Del \vtB$ belongs to 
$L^q \left( \Omega_h \times [-1,1]_{\{ -1, 1\} } \right)$ for any $1 \leq q < \infty$. 
Accordingly, we conclude 
\begin{equation} \label{a4}
\tvt_B \in W^{2,q}(\Omega) \ \mbox{for any} \ 1 \leq q < \infty,
\end{equation}	
in particular, 
\begin{equation} \label{a5}
\| \Grad \widehat{\vt}_B \|_{L^\infty(\Omega; R^3)} \aleq 1.	
\end{equation}
Here and hereafter, we use the symbol $a \aleq b$ to say there is a positive constant $C$ such that $a \leq Cb$.

We will not distinguish between the boundary 
data $\vtB$ and their harmonic extension $\widehat{\vt}_B$ denoting both by the same symbol $\vtB$.

\subsubsection{Energy bounds}
\label{EB}

To eliminate the problem of non--local boundary conditions, we introduce a new quantity 
\begin{equation} \label{a5a}
Z = \Theta + \alpha  \avintO{ \Theta } - \vt_B 
\end{equation}
Computing 
\[
\partial_t Z = \partial_t \Theta + \alpha \avintO{ 
	\partial_t \Theta } 
\ \Rightarrow 
\partial_t \Theta = \partial_t \left( Z - \frac{\alpha}{1 + \alpha} 
 \avintO{ Z } \right)
\]
we may rewrite equation \eqref{t3} in the form
\begin{align} 
\partial_t \left( Z - \frac{\alpha}{\alpha + 1} 
\avintO{ Z } \right) + \Div (\vuh (Z + \vt_B)) &= 
\kappa \Del Z \ \mbox{in}\ (0,T) \times \Omega, \br
Z|_{\partial \Omega} = 0.
\label{a5b}
	 \end{align}
	 
Multiplying \eqref{a5} on $Z$ and integrating over $\Omega$, we deduce the standard ``thermal'' energy balance
\begin{align} 
\frac{1}{2} &\frac{\D }{\dt} \left[ \avintO{ Z^2 } - 
\frac{\alpha}{1 + \alpha} \left( \avintO{ Z } \right)^2  \right]+ {\kappa}\avintO{ |\Grad Z |^2 } \br
&= \avintO{ \vt_B \vuh \cdot \Grad Z }. 
\label{a6}
\end{align}	
	 
Similarly, the scalar product with $\vuh$ and integration by parts in the momentum equation
\eqref{t2} yields
\begin{equation} \label{a7}
\frac{\D }{\dt} \int_{\Omega_h} \frac{1}{2} |\vuh |^2 \ \D x_h + 
\mu \int_{\Omega_h} |\Gradh \vuh |^2 \ \D x_h = \int_{\Omega_h} \left< \Theta \right> \Gradh F \cdot \vuh \ \D x_h	.
\end{equation}

Combining \eqref{a6}, \eqref{a7} we obtain the standard energy {\it a priori} bounds 
\begin{align} 
\sup_{\tau \in [0,T]} \| \vuh (\tau, \cdot) \|_{L^2(\Omega_h; R^2)} &\aleq C\left(\| \vu_{0,h} \|_{L^2(\Omega_h; R^2)}, \| \Theta_0 \|_{L^2(\Omega)}, T, \data
\right) , \br \sup_{\tau \in [0,T]} \| \Theta (\tau, \cdot) \|_{L^2(\Omega)} &\aleq C\left(\| \vu_{0,h} \|_{L^2(\Omega_h; R^2)}, \| \Theta_0 \|_{L^2(\Omega)}, T, \data
\right) , \br
\| \vuh \|_{L^2(0,T; W^{1,2}_0 (\Omega_h; R^2))} & \aleq C\left(\| \vu_{0,h} \|_{L^2(\Omega_h; R^2)}, \| \Theta_0 \|_{L^2(\Omega)}, T, \data
\right) , \br
\| \Theta \|_{L^2(0,T; W^{1,2} (\Omega))} & \aleq C\left(\| \vu_{0,h} \|_{L^2(\Omega_h; R^2)}, \| \Theta_0 \|_{L^2(\Omega)}, T, \data
\right) .
	\label{a8}
\end{align}	
Here and hereafter, ``data'' include exclusively
\begin{equation} \label{a8a}
\| \nabla_h F \|_{L^\infty(\Omega, R^2)},\ \| \vtB \|_{W^{2,q}(\Omega)}, \alpha, a 
\end{equation}
as well as various constants depending on the geometry of $\Omega_h$.

\subsubsection{Maximum principle for the temperature}

In view of \eqref{a8}, the function
\[
t \mapsto \avintO{ \Theta (t,x)  } 
\]
is bounded on compact time intervals. Going back to equation \eqref{t3} 
we may therefore use the standard maximum principle obtaining 
\begin{equation} \label{a10}
\sup_{\tau \in [0,T]} \| \Theta (\tau, \cdot) \|_{L^\infty(\Omega)} 
\leq C\left(\| \vu_{0,h} \|_{L^2(\Omega_h; R^2)}, \| \Theta_0 \|_{L^\infty(\Omega)}, T, \data
\right).
\end{equation}

\subsubsection{Regularity of the velocity field}

In view of the uniform bound \eqref{a10}, the velocity field 
$\vuh$ satisfies the $2-D$ Navier--Stokes system 
\eqref{t2}, with uniformly bounded driving force. Thus we may use the nowadays standard $L^p-$theory (see  Gerhardt \cite{Gerh}, Giga, Miyakawa \cite{GigMiy}, Solonnikov \cite{Solon}) 
to conclude
\begin{align}
\| \partial_t \vuh \|_{L^p(0,T; L^p(\Omega_h; R^2))} & \aleq C\left(\| \vu_{0,h} \|_{W^{2 - \frac{2}{p},p}(\Omega_h; R^2)}, \| \Theta_0 \|_{L^\infty(\Omega)}, T, \data
\right)\br 
\| \vuh \|_{L^p(0,T; W^{2,p}(\Omega_h; R^2))} & \leq 
C\left(\| \vu_{0,h} \|_{W^{2 - \frac{2}{p},p}(\Omega_h; R^2)}, \| \Theta_0 \|_{L^\infty(\Omega)}, T, \data
\right) 
\label{a11}
\end{align}
whenever $2 \leq p < \infty$.

\subsubsection{Regularity of the time derivative of the temperature}

Next, we derive estimates on the time derivative of $\Theta$. 
Denoting $V = \partial_t Z$, where $Z$ has been introduced in 
\eqref{a5a}  we may differentiate \eqref{a5b} obtaining 
\begin{align} 
	\partial_t \left( V - \frac{\alpha}{1 + \alpha} 
	\frac{1}{|\Omega|} \intO{ V } \right) + \Div (\vuh V) &= 
	\kappa \Del V - \Div (\partial_t \vuh (Z + \vt_B )) , \br
	V|_{\partial \Omega} = 0.
	\label{a12}
\end{align}	
Using regularity $\vuh$ established in \eqref{a11}, we obtain the standard ``energy'' estimates 
\begin{align}
V \ \mbox{bounded in}\ L^\infty(0,T; L^2(\Omega; R^3)) \cap 
L^2(0,T; W^{1,2}(\Omega; R^3)), 	
	\nonumber
\end{align}
which yields
\begin{align}
	\partial_t \Theta \ \mbox{bounded in}\ L^\infty(0,T; L^2(\Omega; R^3)) \cap 
	L^2(0,T; W^{1,2}(\Omega; R^3)). 	
	\nonumber
\end{align}
Of course, to perform this step we need $\partial_t \Theta (0, \cdot)$ square integrable. Accordingly the estimates depend on 
$\| \Theta_0 \|_{W^{2,2}(\Omega)}$: 
\begin{align}
\sup_{\tau  \in [0,T]} \| \partial_t \Theta(\tau, \cdot) \|_{L^2(\Omega)} &\leq C\left(\| \vu_{0,h} \|_{W^{2 - \frac{2}{p},p}(\Omega_h; R^2)}, \| \Theta_0 \|_{W^{2,2}(\Omega)}, T, \data
\right)	 \br 
 \| \Theta \|_{L^2(0,T; W^{1,2}(\Omega))} &\leq 
 C\left(\| \vu_{0,h} \|_{W^{2 - \frac{2}{p},p}(\Omega_h; R^2)}, \| \Theta_0 \|_{W^{2,2}(\Omega)}, T, \data
\right). 
\label{a13}	
\end{align}	

\subsubsection{$L^p$-regularity of the temperature}

Going back to \eqref{a5b} we deduce 
\begin{align} \label{a14}
\partial_t Z - \kappa \Del Z &= - \Div (\vuh Z) + \xi, \br 
Z|_{\partial \Omega} &= 0,  	
	\end{align}
where, 
\[
Z = \Theta + \alpha \frac{1}{|\Omega|} \intO{ \Theta } - \vt_B. 
\]
Moreover, 
in view of the previous estimates 
\begin{equation} \label{a15}
	Z \in L^q((0,T) \times \Omega) \ \mbox{for any}\ 1 \leq q < \infty. 
\end{equation}

Now, similarly to Section \ref{h}, we may extend the right--hand side of 
\eqref{a14} as an odd function on the domain 
$\Omega_h \times [-1,1]|_{\{ -1; 1\} }$ to deduce successively the estimates 
\begin{align} 
	\| \partial_t \Theta \|_{L^p(0,T; L^p(\Omega))} & \aleq C\left(\| \vu_{0,h} \|_{W^{2 - \frac{2}{p},p}(\Omega_h; R^2)}, \| \Theta_0 \|_{W^{2 - \frac{2}{p},p}(\Omega)}, T, \data
	\right),\br 
	\| \Theta \|_{L^p(0,T; W^{2,p}(\Omega))} & \aleq C\left(\| \vu_{0,h} \|_{W^{2 - \frac{2}{p},p}(\Omega_h; R^2)}, \| \Theta_0 \|_{W^{2 - \frac{2}{p},p}(\Omega; R^2)(\Omega)}, T, \data
	\right),\ 
	2 \leq p < \infty.
	\label{a16}
\end{align}

We have shown all {\it a priori} estimates compatible with regularity 
of strong solutions claimed in Theorem \ref{eP1}. The solutions can be constructed by means of the mixed Faedo-Galerkin method 
exactly as in \cite{AbbFei22}. We have proved Theorem \ref{eP1}.

\section{Long time behaviour}
\label{L}

We show the rOB system is dissipative in the sense of Levinson, meaning there is a bounded absorbing set the size of which is independent of the initial data. 
To this end, we make an extra assumption 
\begin{equation} \label{L1}
	0 < \alpha < 1 
\end{equation}
in the boundary condition \eqref{i8}. Since 
\[
\alpha = \frac{c_p(\Ov{\vr}, \Ov{\vt} )}{c_v(\Ov{\vr}, \Ov{\vt} )} - 1
\] 
(see \cite{FanFei2025}), hypothesis \eqref{L1} is physically admissible at least for gases.

\subsection{Maximum principle revisited}

Applying the standard maximum principle to the parabolic equation \eqref{t3} we deduce that 
the maximum of $\Theta$ in $[0,T] \times \Ov{\Omega}$ is either attained 
at $t = 0$, in which case
\begin{equation} \label{L2}
\max_{t \in [0,T], x \in \Ov{\Omega}} \Theta(t,x) = \max_{x \in \Ov{\Omega}} \Theta_0,	
\end{equation}	
or there exist $0 < \tau \leq T$, $\Ov{x} \in \partial \Omega$ 
such that
\[
\max_{t \in [0,T], x \in \Ov{\Omega}} \Theta(t,x) = \Theta(\tau, \Ov{x}) = \vtB(\Ov{x}) - \alpha \avintO{ \Theta (\tau, \cdot) }
\leq \max_{x \partial \Omega} \vtB(x) + \alpha \max_{t \in [0,T], x \in \Ov{\Omega}} \Theta(t,x).
\]
As $\alpha < 1$ and $T$ arbitrary we obtain 
\[ 
\max_{t \geq 0,\ x \in \Ov{\Omega}}	\Theta(t,x) \leq 
\max \left\{ \| \Theta_0 \|_{L^\infty(\Omega)}; \frac{1}{1 - \alpha} \| \vtB \|_{L^\infty(\partial \Omega)} \right\}.
\]
	
Repeating the same argument for minimum, we conclude
\begin{equation} \label{L3}
\| \Theta \|_{L^\infty((\tau, \infty) \times \Omega)} \leq \max \left\{ \| \Theta(\tau, \cdot) \|_{L^\infty(\Omega)}; \frac{1}{1 - \alpha} \| \vtB \|_{L^\infty(\partial \Omega)} \right\}.
\end{equation}	
In contrast with \eqref{a10}, estimate \eqref{L3} is global in time.

\subsection{Bounded absorbing set in $L^2 \times L^2$}

First we show that solutions become ultimately bounded in $L^2$. 

\begin{Proposition}[{\bf Absorbing set in $L^2 \times L^2$}] \label{LP1}
Let 
\[
0 < \alpha < 1. 
\]	
Under the hypotheses of Theorem \ref{eP1}, there exists a universal constant $\mathcal{E}_{2,2}$ such that 
\begin{equation} \label{L4}
\limsup_{t \to \infty} \Big[ \| \vuh(t, \cdot) \|_{L^2(\Omega_h; R^2)} + 
\| \Theta(t, \cdot) \|_{L^2(\Omega)} \Big] \leq \mathcal{E}_{2,2}  			
\end{equation}	
for any global--in--time solution of the rOB system. 
The constant $\mathcal{E}_{2,2}$ depends solely on the data as specified in \eqref{a8a}. In particular, it is independent of the initial data.
	
\end{Proposition}

The proof of Proposition \ref{LP1} is based on the property of Levinson dissipativity that can be shown exactly as in \cite[Theorem 4.1]{FeRoSc2024}.

\begin{Proposition}[{\bf Levinson dissipativity}] \label{LP2}	
	Suppose 
\[
0 <\alpha < 1. 
\] 
Let 
$(\vu_{h,n}, \Theta_n)_{n=1}^\infty$ be a sequence of solutions to system \eqref{t2}--\eqref{t4}
defined on the time intervals $(T_n, \infty)$, $T_n \to - \infty$ such that 
\[
	\| \vu_{h,n}(T_n, ,\cdot) \|_{L^2(\Omega_h; R^2)} +  \| \Theta_n(T_n, \cdot) \|_{L^\infty(\Omega)} \leq \mathcal{F}_0\ \mbox{uniformly for}\ n \to \infty.	
\]	

Then, up to a suitable subsequence, 
\begin{align}
	\vun &\to \vuh \ \mbox{in}\ C_{\rm weak, loc}(R; L^2(\Omega; R^d)) \ \mbox{and weakly in}\ L^2_{\rm loc}(R; W^{1,2}_0 (\Omega; R^d)), \br 
	\Theta_n &\to \Theta \ \mbox{in}\ C_{\rm loc}(R; L^2(\Omega)),
	\ \mbox{weakly-(*) in}\ L^\infty_{\rm loc}(R \times \Ov{\Omega}), 
	\ \mbox{and weakly in} \ L^2_{\rm loc}(R; W^{1,2} (\Omega)), 
	\label{L5}
\end{align}
where $(\vuh, \Theta)$ is a (weak) solution of system \eqref{t2}--\eqref{t4} defined in $R \times \Omega$ and satisfying
\begin{equation} \label{L6}
	\| \vuh(t, \cdot) \|_{L^2(\Omega_h; R^2)} + \| \Theta (t, \cdot) \|_{L^\infty(\Omega)} \leq \mathcal{E}_{2,\infty},
\end{equation}
where $\mathcal{E}_{2,\infty}$ depends solely on the data.  In particular, it is independent of $\mathcal{F}_0$.
\end{Proposition}

\begin{Remark} \label{LR1}
	In comparison with Proposition \ref{LP1}, the hypotheses of Proposition \eqref{LP2} include the boundedness 
	of initial data at $T_n$ by the constant $\mathcal{F}_0$. Thus, as observed in \cite[Section 4]{FeRoSc2024},   
	Proposition \ref{LP2} is not a statement about bounded absorbing set in $L^2 \times L^\infty$ but rather
	about regularity of the global attractor. 
\end{Remark}

As shown in \cite[Corollary 4.3]{FeRoSc2024}, Proposition \ref{LP2} implies Proposition \ref{LP1}.

\subsection{Ultimate boundedness of the velocity} 

The uniform bounds claimed in Proposition \ref{LP1}, together with the energy balances \eqref{a6}, \eqref{a7} yield 
\begin{align} 
	\frac{1}{2} &\frac{\D }{\dt} \left[ \avintO{ Z^2 } - 
	\frac{\alpha}{1 + \alpha} \left( \avintO{ Z } \right)^2  \right]+ {\kappa}\avintO{ |\Grad Z |^2 } \br
	&+ \frac{\D }{\dt} \int_{\Omega_h} \frac{1}{2} |\vuh |^2 \ \D x_h + 
	\mu \int_{\Omega_h} |\Gradh \vuh |^2 \ \D x_h \leq K(t), 
	\label{L7}
\end{align}	
where 
\begin{equation} \label{L8}
\limsup_{t \to \infty} K(t) \leq K_\infty,\ K_\infty = K_\infty(\data).	
\end{equation}	
From this and Proposition \ref{LP1}, we deduce 
\begin{equation} \label{L9}
\limsup_{T \to \infty} \left[ \int_T^{T+1} \| \vuh \|^2_{W^{1,2}(\Omega_h; R^2)} \dt +  
\int_T^{T+1} \| \Theta \|^2_{W^{1,2}(\Omega)} \dt \right] \leq \mathcal{E}_{1,2;1,2}
\end{equation}

\subsubsection{Maximal regularity estimates for the Navier--Stokes system in $2D$}

We make use of the exact formula derived by Gerhardt \cite[Theorem 3]{Gerh}:
\begin{align} 
&\int_{T}^{T+1} \left[ \| \vuh \|^2_{W^{2,2}(\Omega_h; R^2)}  + \| \partial_t \vuh \|^2_{L^{2}(\Omega_h; R^2)}\right] \dt \br 
&\leq c_1 \left( \| \vuh(T) \|_{W^{1,2}(\Omega_h; R^2)} + \int_T^{T+1} \int_{\Omega_h} |\left< \Theta \right>|^2 |\Gradh \vc{G}|^2 \dx_h \dt \right) \times \br 
&\times \left[ 1 + c_0 c_1 \left( \sup_{t \in [T, T+1]} \| \vuh \|^2_{L^2(\Omega_h; R^2)} \int_T^{T+1} \| \vuh \|^2_{W^{1,2}(\Omega_h; R^2)}   \right) \times 
\right. \br 
&\times \left. 
\exp \left( c_0 c_1^2 \left( \sup_{t \in [T, T+1]} \| \vuh \|^2_{L^2(\Omega_h; R^2)} \int_T^{T+1} \| \vuh \|^2_{W^{1,2}(\Omega_h; R^2)}   \right)    \right) \right],
\label{L10}
\end{align}	
where $c_0$, $c_1$ depend only on the domain $\Omega_h$.
This estimate, together with \eqref{L4}, \eqref{L9} yields 
\begin{equation} \label{L12}
\limsup_{T \to \infty} \left[ \| \vuh(T, \cdot) \|^2_{W^{1,2}(\Omega_h; R^2)} +\int_{T}^{T+1} \left[ \| \vuh \|^2_{W^{2,2}(\Omega_h; R^2)}  + \| \partial_t \vuh \|^2_{L^{2}(\Omega_h; R^2)}\right] \dt \right] \leq \mathcal{E}_{u,1,2}.
\end{equation}

\subsection{Ultimate boundedness of the temperature deviation}

Finally, using the estimates \eqref{L12} and multiplying \eqref{a14} by $\Del Z$, we deduce 
\begin{equation} \label{L13}
\limsup_{T \to \infty} \int_T^{T+1} \| \Theta \|^2_{W^{2,2}(\Omega)} \leq \mathcal{E}_{2,2}.
\end{equation}	
As $W^{2,2} \hookrightarrow L^\infty$, the estimate \eqref{L13} along with the maximum principle yield 
\begin{equation} \label{L14}
\limsup_{t \to \infty} \| \Theta (t, \cdot) \|_{L^\infty(\Omega)} \leq \mathcal{E}_{\infty}. 	 
\end{equation}	

Let us summarize the results obtained in this section.

\begin{mdframed}[style=MyFrame]

\begin{Theorem}[{\bf Absorbing set in $W^{1,2} \times L^\infty$}] \label{LP3}
Let 
\[
0 < \alpha < 1. 
\]	
Under the hypotheses of Theorem \ref{eP1}, there exists a universal constant $\mathcal{E}_{1,2; \infty}$ such that 
\begin{equation} \label{L15}
	\limsup_{t \to \infty} \Big[ \| \vuh(t, \cdot) \|_{W^{1,2}_0(\Omega_h; R^2)} + 
	\| \Theta(t, \cdot) \|_{L^\infty \cap W^{1,2}(\Omega)} \Big] \leq \mathcal{E}_{1,2; \infty}  			
\end{equation}	
for any global--in--time solution emanating from the initial data satisfying the hypotheses of Theorem \ref{eP1}. 
The constant $\mathcal{E}_{1,2; \infty}$ depends solely on the data as specified in \eqref{a8a}. In particular, it is independent of the initial data.

\end{Theorem}

\end{mdframed}

Using Theorem \ref{LP3}, together with Proposition \ref{LP2}, and the uniform bounds established in \eqref{L12}, \eqref{L13}, we deduce, 
exactly as in \cite{FeRoSc2024}, that the rOB system admits a global trajectory attractor 
\begin{align}
\mathcal{A} = &\left\{ (\vuh, \Theta) \in BC \Big(R; W^{1,2}_0 (\Omega_h; R^2) \times (W^{1,2} \cap L^\infty) (\Omega) \Big) \ \Big| \right. \br &(\vuh, \Theta) 
\ \mbox{solves the rOB system in} \ R \times (\Omega_h \times \Omega)
\br
&\sup_{\tau \in R} \left[ \| \vuh(\tau, \cdot) \|_{W^{1,2}(\Omega_h; R^2)} + \| \Theta(\tau, \cdot) \|_{W^{1,2}(\Omega)} +
\| \Theta(\tau, \cdot) \|_{L^\infty(\Omega)} \right] \leq \mathcal{E}({\rm data}), \br 
&\sup_{\tau \in R} \int_{\tau}^{\tau + 1} \left[ \|  \partial_t \vuh \|_{L^{2}(\Omega_h; R^2)} + \|  \vuh \|_{W^{2,2}(\Omega_h; R^2)} 
+ \|  \partial_t \Theta \|_{L^{2}(\Omega_h)} + \|  \Theta \|_{W^{2,2}(\Omega)}\right] \dt \leq \mathcal{E}({\rm data}) \Big\}.
\label{AT}
\end{align}
The attractor $\mathcal{A}$ is non--empty and time-shift invariant.

\section{Continuous data assimilation method for rotating Ober\-beck--Boussinesq  system}
\label{D}

Following the idea of Azouani, Olson, and Titi \cite{AzOlTi}, see also Carlson, Larios, and Titi \cite{CaLaTi}, Farhat, Jolly, M. S., and Titi \cite{FaJoTi}, 
we develop a data assimilation method for the rOB system. 

We start by introducing the interpolant operators
\begin{align} 
I_\delta [\vc{v}] ,\ \| I_\delta[\vc{v}] \|_{L^2(\Omega_h; R^2)} &\leq \Ov{I} \Big(\| \vc{v} \|_{L^2(\Omega_h;R^2)} \Big),\ 
\| I_\delta[\vc{v}] - \vc{v} \|_{L^2(\Omega_h; R^2)} \leq \mathcal{O}(\delta) \| \vc{v} \|_{W^{1,2}(\Omega_h; R^2)},\br 
\mathcal{O}(\delta) &\to 0 \ \mbox{as}\ \delta \to 0.\label{D1}
\end{align}

Next, consider an entire solution of the rOB system $(\vuh, \Theta) \in \mathcal{A}$, meaning 
\begin{align} 
	\Divh \vuh &= 0, \br
	\partial_t \vuh + \Divh( \vuh \otimes \vuh) + \Gradh {\Pi} &=
	\mu \Delta_h \vuh - \left< \Theta \right > \Gradh F
	\label{D2} 
\end{align}
in $R \times \Omega_h$,
\begin{equation} \label{D3} 
	\partial_t \Theta + \Div (\Theta \vuh) = \kappa \Del \Theta 
\end{equation}
in $R \times \Omega$,  
\begin{equation} \label{D4}
	\vuh|_{\partial \Omega_h} = 0,\ 
	\Theta|_{\partial \Omega} = \vtB - \alpha \avintO{ \Theta }. 
\end{equation}

Now, we introduce a ``forced'' rOB system 
\begin{align} 
	\Divh \tvuh &= 0, \br
	\partial_t \tvuh + \Divh( \tvuh \otimes \tvuh) + \Gradh \widetilde{\Pi }&=
	\mu \Delta_h \tvuh - \Lambda I_\delta[ \tvuh - \vuh ] - \left< \tTheta \right > \Gradh G
	\label{D5} 
\end{align}
in $(0,T) \times \Omega_h$,
\begin{equation} \label{D6} 
	\partial_t \tTheta + \Div (\tTheta \tvuh) = \kappa \Del \tTheta 
\end{equation}
in $(0,T) \times \Omega$, 
with the boundary conditions 
\begin{equation} \label{D7}
	\tvuh|_{\partial \Omega_h} = 0,\ 
	\tTheta|_{\partial \Omega} = \vtB - \alpha \avintO{ \tTheta },  
\end{equation}
and the initial data 
\begin{equation} \label{D8}
\tvuh(0, \cdot) = \tvu_{0,h},\ \tTheta(0, \cdot) = \tTheta_0.	
\end{equation}	
In view of \eqref{D1}, global solvability of the forced problem \eqref{D5}--\eqref{D8} can be shown repeating the arguments of the proof 
of Theorem \ref{eP1}. Here $\Lambda > 0$ is the so-called nudging parameter to be determined below.

Our ultimate task is to compare the distance of the solutions $(\tvuh, \tTheta)$, $(\vuh, \Theta)$.
Similarly to Section \ref{EB}, we introduce 
\[
	Z = \Theta + \alpha  \avintO{ \Theta } - \tvt_B,\ \tvZ = \tTheta + \alpha  \avintO{ \tTheta }  - \vt_B
\]
rewriting \eqref{D3}, \eqref{D6} in the form 
\begin{align} 
	\partial_t \left( \tvZ - \frac{\alpha}{\alpha + 1} 
	\avintO{ \tvZ } \right) + \Div (\tvuh ( \tvZ + \vt_B)) &= 
	\kappa \Del \tvZ \ \mbox{in}\ R \times \Omega, \br
	\tvZ|_{\partial \Omega} = 0, 
	\label{D9}
\end{align}
\begin{align} 
	\partial_t \left( Z - \frac{\alpha}{\alpha + 1} 
	\avintO{ Z } \right) + \Div (\vuh ( Z + \vt_B)) &= 
	\kappa \Del Z \ \mbox{in}\ (0,T) \times \Omega, \br
	Z|_{\partial \Omega} = 0, 
	\label{D10}
\end{align}

Subtracting the corresponding equations we get 
\begin{align} 
	\Divh (\tvuh - \vuh) &= 0, \br
	\partial_t (\tvuh - \vuh)  + \Gradh P &=
	\mu \Delta_h (\tvuh - \vuh) - \Lambda I_\delta[ \tvuh - \vuh ] - \left< \tvZ - Z \right > \Gradh G\br &+ \Divh( \vuh \otimes \vuh)-  \Divh( \tvuh \otimes \tvuh), \br (\tvuh - \vuh)|_{\partial \Omega_h} &= 0,
	\label{D11} 
\end{align}
and, 
\begin{align} 
	\partial_t \left( (\tvZ-Z) - \frac{\alpha}{\alpha + 1} 
	\avintO{ (\tvZ-Z) } \right) + \Div (\tvuh ( \tvZ-Z)) &= 
	\kappa \Del (\tvZ-Z) + \Div ( \tvZ (\tvuh - \vuh) ) , \br
	(\tvZ-Z)|_{\partial \Omega} &= 0. 
	\label{D12}
\end{align}

Multiplying equation \eqref{D11} on $(\tvuh - \vuh)$ and integrating by parts we get 
\begin{align} 
\frac{1}{2} \frac{\D}{\dt} \| \tvuh - \vuh \|^2_{L^2(\Omega_h; R^2)}  &+ \mu \| \nabla_h( \tvuh - \vuh) \|^2_{L^2(\Omega_h; R^4)}
+ \Lambda \intOh{ I_\delta[ \tvuh - \vuh ] \cdot (\tvuh - \vuh) }\br = &\intOh{ \left< Z-\tvZ \right > \Gradh G \cdot (\tvuh - \vuh) } \br
&+ \intOh{ \Divh( \vuh \otimes \vuh - \tvuh \otimes \tvuh) \cdot (\tvuh - \vuh) }
\label{D13}
\end{align}	
Similarly, multiplying \eqref{D12} by $(\tvZ-Z)$ and integrating by parts yields 
\begin{align} 
	\frac{1}{2} &\frac{\D }{\dt} \left[ \avintO{ (\tvZ-Z)^2 } - 
	\frac{\alpha}{1 + \alpha} \left( \avintO{ (\tvZ-Z) } \right)^2  \right]+ {\kappa}\avintO{ |\Grad (\tvZ-Z) |^2 } \br
	&= \avintO{ Z (\tvuh - \vuh) \cdot \Grad (\tvZ-Z)  }. 
	\label{D14}
\end{align}	

In accordance with \eqref{D1}, we have 
\begin{align}
	&\intOh{ I_\delta[ \tvuh - \vuh ] \cdot (\tvuh - \vuh) } = \intOh{ ( \tvuh - \vuh )^2 } \br
	&+ \intOh{ \Big( I_\delta[ \tvuh - \vuh ] - (\tvuh - \vuh) \Big) \cdot (\tvuh - \vuh) }, \br
	&\left| \intOh{ \Big( I_\delta[ \tvuh - \vuh ] - (\tvuh - \vuh) \Big) \cdot (\tvuh - \vuh) } \right| \leq \mathcal{O}(\delta) 
	\| \tvuh - \vuh \|_{W^{1,2}(\Omega_h; R^2)} \| \| \tvuh - \vuh \|_{L^{2}(\Omega_h; R^2)}. 
\label{D15}
\end{align}	
Moreover, we rewrite 
\begin{align}
	\intOh{ \Divh( \vuh \otimes \vuh - \tvuh \otimes \tvuh) \cdot (\tvuh - \vuh) } = \intOh{ \nabla_h \vuh : (\tvuh - \vuh) \otimes (\tvuh - \vuh) },
\nonumber
\end{align}	
where, by Ladyzhenskaya interpolation inequality 
\begin{align} 
&\left| \intOh{ \nabla_h \vuh : (\tvuh - \vuh) \otimes (\tvuh - \vuh) } \right| 
\leq \| \nabla_h \vuh \|_{L^2(\Omega_h; R^4)} \| \tvuh - \vuh \|^2_{L^4(\Omega_h; R^2)}\br &\leq C_L  \| \nabla_h \vuh \|_{L^2(\Omega_h; R^4)} 
 \| \tvuh - \vuh \|_{L^2(\Omega_h; R^2)} \| \nabla_h (\tvuh - \vuh) \|_{L^2(\Omega_h; R^4)} .
\label{D16}
\end{align}
Finally, we recall Poincar\' e inequality 
\[
\| \tvuh - \vuh \|^2_{W^{1,2}_0(\Omega_h; R^2)} \leq C_P \intOh{ |\nabla_h (\tvuh - \vuh) |^2 },\ 
\| \tvZ-Z \|^2_{W^{1,2}_0(\Omega)} \leq C_P \intOh{ |\Grad (\tvZ-Z) |^2 }. 
\]

Summing up identities \eqref{D13}, \eqref{D14} and making use of the estimates \eqref{D15}, \eqref{D16},  we obtain
\begin{align} 
\frac{1}{2} &\frac{\D }{\dt} \left(  \| \tvuh - \vuh \|^2_{L^2(\Omega_h; R^2)} + \left[ \avintO{ (\tvZ-Z)^2 } - 
\frac{\alpha}{1 + \alpha} \left( \avintO{ (\tvZ-Z) } \right)^2  \right] \right) \br 
&+ \mu \| \nabla_h( \tvuh - \vuh) \|^2_{L^2(\Omega_h; R^4)}+ \Lambda  \|  \tvuh - \vuh \|^2_{L^2(\Omega_h; R^2)} + {\kappa}\avintO{ |\Grad (\tvZ-Z) |^2 } \br
&\leq \mathcal{O}(\delta) \Lambda
\| \tvuh - \vuh \|_{W^{1,2}(\Omega_h; R^2)} \  \| \tvuh - \vuh \|_{L^{2}(\Omega_h; R^2)} \br 
&+ C_L  \| \nabla_h \vuh \|_{L^2(\Omega_h; R^4)} 
\| \tvuh - \vuh \|_{L^2(\Omega_h; R^2)} \| \nabla_h (\tvuh - \vuh) \|_{L^2(\Omega_h; R^4)} \br 
&+ \| Z \|_{L^\infty(\Omega)} \| \tvuh - \vuh \|_{L^{2}(\Omega_h; R^2)} \| \Grad (\tvZ-Z) \|_{L^{2}(\Omega; R^3)} \br 
&+ \| \nabla_h F \|_{L^\infty(\Omega_h; R^2)} \| \tvuh - \vuh \|_{L^{2}(\Omega_h; R^2)} \| \tvZ-Z \|_{L^{2}(\Omega)} 
\label{D17}
\end{align}	
Thus we can choose first $\Lambda = \Lambda( \mathcal{E}, C_L, C_P, \data )$ large and 
$\delta = \delta(\Lambda) > 0$ small enough so that inequality \eqref{D17} reduces to 
\begin{align} 
	\frac{\D }{\dt} &\left(  \| \tvuh - \vuh \|^2_{L^2(\Omega_h; R^2)} + \left[ \avintO{ (\tvZ-Z)^2 } - 
	\frac{\alpha}{1 + \alpha} \left( \avintO{ (\tvZ-Z) } \right)^2  \right] \right) \br 
	&+ \beta \left[ \|  \tvuh - \vuh \|^2_{W^{1,2}_0(\Omega_h; R^2)}  + \|  \tvZ-Z \|^2_{W^{1,2}_0(\Omega)} \right]\leq 0, 
	\label{D18}
\end{align}	
where $\beta = \beta(\delta, \Lambda, \mathcal{E}, C_L, C_P, \data) > 0$.

Note that 
\begin{align}
 \avintO{ (\Theta - \tTheta)^2}&=
 \avintO{ (\tvZ-Z)^2 } - \frac{2\alpha}{1 + \alpha} \left( \avintO{ (\tvZ-Z) } \right)^2  
+\left(\frac{\alpha}{1 + \alpha} \right)^2\left( \avintO{ (\tvZ-Z) } \right)^2  \\
&\le \avintO{ (\tvZ-Z)^2 } - 
	\frac{\alpha}{1 + \alpha} \left( \avintO{ (\tvZ-Z) } \right)^2.\nonumber
\end{align}
In addition, Jensen's inequality gives
\begin{align}
 \avintO{ (\tvZ-Z)^2 } - 
	\frac{\alpha}{1 + \alpha} \left( \avintO{ (\tvZ-Z) } \right)^2
\le \frac{1}{1+\alpha}\avintO{ (\tvZ-Z)^2 }
\end{align}

We have shown the following result. 

\begin{mdframed}[style=MyFrame]
	
\begin{Theorem}[{\bf Convergence of data assimilation method}] \label{tt1}
	
Let $(\vuh, \Theta) \in \mathcal{A}$ be an entire solution and let 
$(\tvuh, \tTheta)$ be a solution of the forced  rOB
system \eqref{D5}--\eqref{D8} emanating from the initial data 
$(\widetilde{\vc{u}}_{0,h}, \tTheta_0)$. 

Then there exist $\Lambda = \Lambda (\mathcal{E}, C_L, C_P, \data) > 0$ large, $\delta = \delta(\Lambda) > 0$ small 
such that 
\begin{align}
\| (\vuh - \tvuh)(\tau, \cdot) \|_{L^2(\Omega_h; R^2)} &+ 
\| (\Theta - \tTheta)(\tau, \cdot) \|_{L^2(\Omega)}\br &\aleq 
\left( 1 + \| \widetilde{\vc{u}}_{0,h} \|_{L^2(\Omega_h; R^2)} +  \| \tTheta_0 \|_{L^2(\Omega)} \right) \exp (- \beta \tau),\ \tau > 0,  
\nonumber
\end{align}
for a certain $\beta = \beta (\mathcal{E}, C_L, C_P, \data)  > 0$. 	
	
\end{Theorem}		
	
\end{mdframed}

\section{Data assimilation via the primitive compressible model}\label{S5}

We suppose that the values of the interpolation operators used in the
forced rOB system \eqref{D5}, \eqref{D6} are provided by the primitive 
system \eqref{s1}--\eqref{s5} in the regime $\ep \to 0$. Accordingly, we replace the forced rOB system \eqref{D5}--\eqref{D7} by 
\begin{align} 
	\Divh \tvuh &= 0, \br
	\Ov{\vr} \partial_t \tvu_h + \Ov{\vr} \Divh( \tvuh \otimes \tvuh) + \Gradh \Pi &=
	\mu \Delta_h \tvuh - 
	\Lambda I_\delta[\tvuh - (\vu)_h]- \left< \widetilde{\theta} \right > \Gradh F
	\label{Di3} 
\end{align}
in $(0,T) \times \Omega_h$,
\begin{equation} \label{Di6} 
	\partial_t \widetilde{\theta} + \Div (\widetilde{\theta} \tvuh) = \kappa \Del \widetilde{\theta} +
	a \Div ({F} \tvuh), \ \kappa > 0,\ a \in R, 
\end{equation}
in $(0,T) \times \Omega$,  with the boundary conditions
\begin{equation} \label{Di8}
\tvuh|_{\partial \Omega_h} = 0,\ 	
	\widetilde{\theta}|_{\partial \Omega} = \vtB - \alpha \frac{1}{|\Omega|} \intO{ 
		\widetilde{\theta} }.
\end{equation}
Here, the symbol $(\vu)_h = (u^1,u^2)$ denotes the horizontal component 
of the velocity $\vu = (u^1,u^2, u^3)$ determined by the primitive system \eqref{s1} -- \eqref{s3}.

\subsection{Initial data for the primitive system} 

The initial data satisfied by the primitive system that determine the 
velocity field $(\vu)_h$ in \eqref{Di3} are {\it a priori} not known. 
As we consider the primitive system in the regime of small $\ep$, we suppose the initial data are \emph{well--prepared}, meaning they do not support acoustic waves. Such an assumption is quite natural for models of atmosphere in meteorology. Accordingly, we suppose there is 
a reference density $\Ov{\vr} > 0$ and a reference temperature 
$\Ov{\vt}$ such that the solution $(\vr, \vt, \vu)$ of the primitive system satisfies
\begin{align}
	\left\| \frac{\vr( 0, \cdot) - \Ov{\vr} }{\ep} \right\|_{L^\infty(\Omega)} \aleq 1 ,\  
	\left\| \frac{\vr( 0, \cdot) - \Ov{\vr}}{\ep} -  {\mathcal{R}}_0 \right\|_{L^1(\Omega)} \leq \ep, \br 
	\left\| \frac{\vt( 0, \cdot) - \Ov{\vt}}{\ep} \right\|_{L^\infty(\Omega)} \aleq 1 ,\  
	\left\| \frac{\vt( 0, \cdot) - \Ov{\vt}}{\ep} -  \mathcal{T}_0 \right\|_{L^1(\Omega)} \leq \ep, \br
	\left\| \vu( 0, \cdot) \right\|_{L^\infty(\Omega; R^3)} \aleq 1,\ 
	\left\| \vu( 0, \cdot) - {\vu}_{0,h} \right\|_{L^1(\Omega; R^3)} \leq \ep
	\label{s12}
\end{align}
for certain functions $
{\mathcal{R}}_0$, $\mathcal{T}_0$, and $\vu_{0,h}$ satisfying
\begin{align} 
	\mathcal{T}_0 &\in W^{2,\infty}(\Omega), \ {\mathcal{T}}_0|_{\partial \Omega} = \vtB, \br 
	[{\vu}_{0,h}, 0]\equiv {\vu}_{0,h} 
	&\in W^{2,\infty}(B(r); R^2),\  
	{\vu}_{0,h}|_{\partial B(r)} = 0,\ \Divh {\vu}_{0,h} = 0, 	
	\label{s13}
\end{align}
and	
\begin{equation} \label{s14}
	\frac{\partial p(\Ov{\vr}, \Ov{\vt})}{\partial \vr } \Grad {\mathcal{R}}_0 + 
	\frac{\partial p(\Ov{\vr}, \Ov{\vt})}{\partial \vt } \Grad {\mathcal{T}}_0 = \Ov{\vr} \Grad F .	
\end{equation}

\subsection{Approximate data assimilation}

Under the hypotheses on the initial data \eqref{s12} -- \eqref{s13}, 
the solutions of the primitive system \eqref{s1}--\eqref{s3} approach the solution of the \eqref{i3}--\eqref{i6} as $\ep \to 0$ 
uniformly on compact time intervals, cf. 
\cite{FanFei2025}. Consequently, combining this result with Theorem \ref{tt1} we obtain the following conclusion. 

\begin{mdframed}[style=MyFrame]
	\begin{Theorem}[{\bf Approximate data assimilation}]
	
Let $\Omega_h = B_r$. Let $(\vr, \vu, \vt)$ be a weak solution of the 
primitive system \eqref{s1}--\eqref{i5} emanating from the initial data 
\eqref{s12}--\eqref{s13}. Let 
$\tvuh$, $\widetilde{\theta}$ by the solution of the forced rOB system 
\eqref{Di3}--\eqref{Di6} emanating from the initial data
$\tvuh(0, \cdot) = \tvu_{0,h}$, $\widetilde{\theta} (0,\cdot) = \widetilde{\theta}_{0}$.  

Then there exists $\Lambda > 0$, $\delta = \delta(\Lambda) > 0$, and 
$\beta > 0$ such that  
\begin{align}
	\| (\vr\vu - \Ov{\vr}\tvuh)(\tau, \cdot) \|_{L^1(\Omega_h; R^2)} &+ 
	\left\| \frac{\vt(\tau, \cdot) - \vtB}{\ep} - \lambda \avintO{\frac{\vt(\tau, \cdot) - \vtB}{\ep} } - \widetilde{\theta} (\tau, \cdot) \right\|_{L^1(\Omega)}\br &\aleq 
	\left( 1 + \| \widetilde{\vc{u}}_{0,h} \|_{L^2(\Omega_h; R^2)} +  \| \widetilde{\theta}_0 \|_{L^2(\Omega)} \right) \exp (- \beta \tau) + 
	\chi(\ep, \tau)
	\nonumber
\end{align}	
where $\chi(\ep, \cdot) \to 0$ as $\ep \to 0$ uniformly on compact time intervals.
\end{Theorem}

\end{mdframed}

\def\cprime{$'$} \def\ocirc#1{\ifmmode\setbox0=\hbox{$#1$}\dimen0=\ht0
	\advance\dimen0 by1pt\rlap{\hbox to\wd0{\hss\raise\dimen0
			\hbox{\hskip.2em$\scriptscriptstyle\circ$}\hss}}#1\else {\accent"17 #1}\fi}


\end{document}